\documentclass[graybox]{svmult}

\usepackage{courier} 
\usepackage[bottom]{footmisc} 
\usepackage{graphicx} 
\usepackage{helvet} 
\usepackage{makeidx} 
\usepackage{mathptmx} 
\usepackage{multicol} 
\makeindex

\begin{document}

\title*{Fast Multipole Method as a Matrix-Free Hierarchical Low-Rank Approximation}
\author{Rio Yokota, Huda Ibeid, David Keyes}
\institute{Rio Yokota \at Tokyo Institute of Technology, 2-12-1 O-okayama Meguro-ku, Tokyo, Japan,
\email{rioyokota@gsic.titech.ac.jp}
\and Huda Ibeid, David Keyes \at King Abdullah University of Science and Technology, 4700 KAUST, Thuwal, Saudi Arabia,
\email{huda.ibeid@kaust.edu.sa,david.keyes@kaust.edu.sa}}
\maketitle

\abstract{There has been a large increase in the amount of work on hierarchical low-rank approximation methods,
where the interest is shared by multiple communities that previously did not intersect.
This objective of this article is two-fold; to provide a thorough review of the recent advancements in this field from both
analytical and algebraic perspectives, and to present a comparative benchmark of
two highly optimized implementations of contrasting methods for some simple yet representative test cases.
We categorize the recent advances in this field from the perspective of compute-memory tradeoff,
which has not been considered in much detail in this area.
Benchmark tests reveal that there is a large difference in the memory consumption and performance
between the different methods.}

\section{Introduction}
The fast multipole method (FMM) was originally developed as an algorithm to bring down the $\mathcal{O}(N^2)$
complexity of the direct $N$-body problem to $\mathcal{O}(N)$
by approximating the hierarchically decomposed far field with multipole/local expansions.
In its original form, the applicability of FMM is limited to problems that have a Green's function solution,
for which the multipole/local expansions can be calculated analytically.
Their function is also limited to matrix-vector multiplications,
in contrast to the algebraic variants that can perform matrix-matrix multiplication and factorizations.
However, these restrictions no longer apply to the FMM
since the kernel independent FMM \cite{Ying2004} does not require a Green's function,
and inverse FMM \cite{Ambikasaran2014} can be used as the inverse operator instead of the forward mat-vec.
Therefore the FMM can be used for a wide range of scientific applications,
which can be broadly classified into elliptic partial differential equations (PDE) and kernel summation.
Integral form of elliptic PDEs can be further categorized into boundary integrals for homogeneous problems,
discrete volume integrals, and continuous volume integrals.

Scientific applications of FMM for boundary integrals include acoustics \cite{Wolf2011a, Hao2015},
biomolecular electrostatics \cite{Yokota2011}, electromagnetics \cite{Darve2004a, Gimbutas2013},
fluid dynamics for Euler \cite{Willis2005} and Stokes \cite{Rahimian2010} flows,
geomechanics \cite{Verde2015}, and seismology \cite{Chaillat2008, Wilkes2015}.
Application areas of FMM for discrete volume integrals are astrophysics \cite{Bedorf2014},
Brownian dynamics \cite{Liang2013}, classical molecular dynamics \cite{Ohno2014},
density functional theory \cite{Shao2001}, vortex dynamics \cite{Yokota2013},
and force directed graph layout\cite{Yunis2012}.
FMM for continuous volume integrals have been used to solve
Schr\"odinger \cite{Zhao2007} and Stokes \cite{Malhotra2014} equations.
More generalized forms of FMM can be used as fast kernel summation for
Bayesian inversion \cite{Ambikasaran2013a}, Kalman filtering \cite{Li2014},
Machine learning \cite{Gray2001,Lee2012}, and radial basis function interpolation \cite{Gumerov2007}.

All of these applications have in common the key feature that they are global problems
where the calculation at every location depends on the values everywhere else.
Elliptic PDEs that represent a state of equilibrium, many iterations with global inner products for their solution,
dense matrices in boundary integral problems, all-to-all interaction in $N$-body problems,
and kernel summations with global support are all different manifestations of the same source of global data dependency.
Due to this global data dependency, their concurrent execution on future computer architectures
with heterogeneous and deep memory hierarchy is one of the main challenges of exascale computing.
For global problems that require uniform resolution, FFT is often the method of choice,
despite its suboptimal communication costs.
The methods we describe here have an advantage for global problems that require non-uniform resolution.
For such non-uniform global problems multigrid methods are known to do quite well.
Whether the reduced synchronization and increased arithmetic intensity of the FMM will become advantageous
compared to multigrid on future architectures is something that is yet to be determined.

Many of the original FMM researchers have now moved on to develop algebraic variants of FMM,
such as $\mathcal{H}$-matrix \cite{Hackbusch1999}, $\mathcal{H}^2$-matrix \cite{Hackbusch2000a},
hierarchically semi-seprable (HSS) \cite{Chandrasekaran2006}, hierarchically block-separable (HBS) \cite{Martinsson2005},
and hierarchically off-diagonal low-rank (HODLR) \cite{Ambikasaran2013} matrices.
The differences between these methods are concisely summarized by Ambikasaran \& Darve \cite{Ambikasaran2014}.
These algebraic generalizations of the FMM can perform addition, multiplication,
and even factorization of dense matrices with near linear complexity.
This transition from analytic to algebraic did not happen suddenly,
and semi-analytic variants were developed along the way \cite{Ying2004, Fong2009}.
Optimization techniques for the FMM such as compressed translation operators and their precomputation,
also fall somewhere between the analytic and algebraic extremes.

The spectrum that spans purely analytic and purely algebraic forms of these hierarchical low-rank approximation methods,
represents the tradeoff between computation (Flops) and memory (Bytes).
The purely analytic FMM is a matrix-free $\mathcal{H}^2$-matrix-vector product,
and due to its matrix-free nature it has very high arithmetic intensity (Flop/Byte) \cite{Barba2013}.
On the other end we have the purely algebraic methods, which precompute and store the entire hierarchical matrix.
This results in more storage and more data movement, both vertically and horizontally in the memory hierarchy.
When the cost of data movement increases faster than arithmetic operations on future architectures,
the methods that compute more to store/move less will become advantageous.
Therefore, it is important to consider the whole spectrum of hierarchical low-rank approximation methods,
and choose the appropriate method for a given pair of application and architecture.

There have been few attempts to quantitatively investigate the tradeoff between the analytic and algebraic
hierarchical low-rank approximation methods.
Previously, the applicability of the analytic variants were limited to problems with Green's functions,
and could only be used for matrix-vector products but not to solve the matrix.
With the advent of the kernel-independent FMM (KIFMM) \cite{Ying2004} and inverse FMM (IFMM) \cite{Ambikasaran2014},
these restrictions no longer apply to the analytic variants.
Furthermore, the common argument for using the algebraic variants because they can operate directly on the matrix
without the need to pass geometric information is not very convincing.
Major libraries like PETSc offer interfaces to insert one’s own matrix free preconditioner as a function,
and passing geometric information is something that users are willing to do if the result is increased performance.
Therefore, there is no strong reason from the user’s perspective to be monolithically inclined to use the algebraic variants.
It is rather a matter of choosing the method with the right balance between its
analytic (Flops) and algebraic (Bytes) features.

The topic of investigating the tradeoff between analytic and algebraic hierarchical low-rank approximation methods
is too broad to cover in a page-constrained article.
In the present work, we limit our investigation to the compute-memory tradeoff
in a comparison between FMM and HSS for Laplace and Helmholtz kernels.
We also investigate the use of FMM as a preconditioner for iterative solutions to
the Laplace and Helmholtz problems with finite elements,
for which we compare with geometric and algebraic multigrid methods.

\begin{figure}[b]
\centering
\includegraphics[scale=.55]{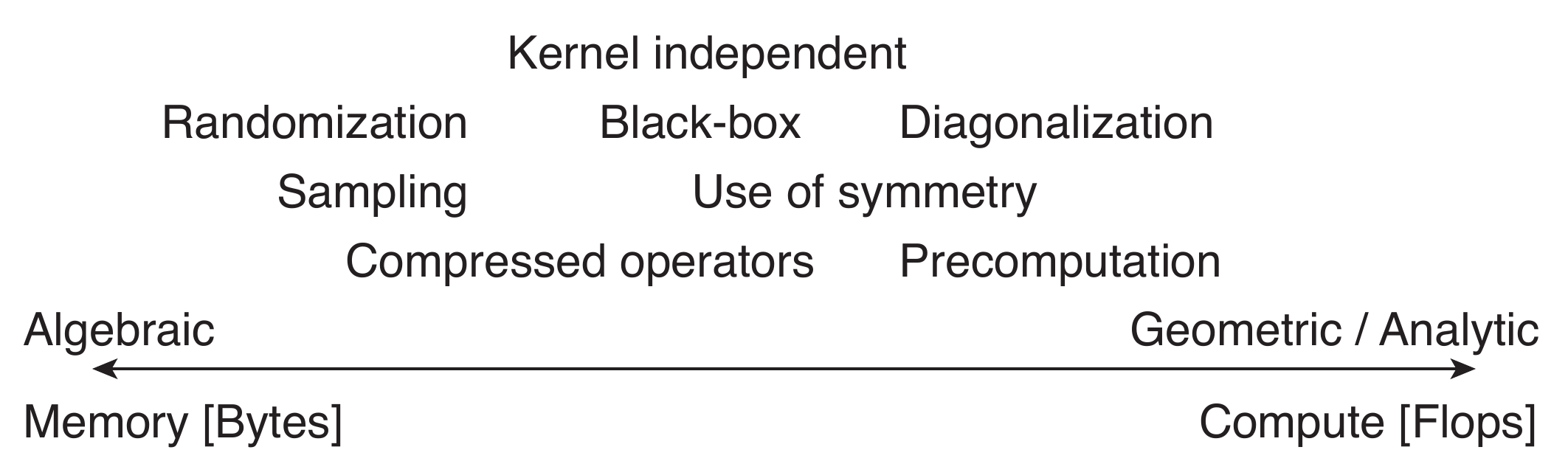}
\caption{The compute-memory tradeoff between the analytic and algebraic hierarchical low-rank approximation methods.
Various techniques lie between the analytic and algebraic extremes.}
\label{fig:tradeoff}
\end{figure}

\section{Hierarchical Low-Rank Approximation: Analytic or Algebraic?}
In this section we review the full spectrum of hierarchical low-rank approximations starting from
the analytic side and proceeding to the algebraic side.
The spectrum is depicted in Fig. \ref{fig:tradeoff},
where various techniques like between the analytic and algebraic extremes. 
One can choose the appropriate method for a given architecture to achieve the best performance.

\subsection{Analytic Low-Rank Approximation}
On the analytic end of the spectrum, we have classical methods such as the Treecode \cite{Barnes1986},
FMM \cite{Appel1985,Greengard1987}, and panel clustering methods \cite{Hackbusch1989}.
These methods have extremely high arithmetic intensity (Flop/Byte) due to their matrix-free nature,
and are compute-bound on most modern architectures.
One important fact is that these are not brute force methods that do unnecessary Flops,
but are (near) linear complexity methods that are only doing useful Flops,
but they are still able to remain compute-bound.
This is very different from achieving high Flops counts on dense matrix-matrix multiplication or LU decomposition
that have $\mathcal{O}(N^3)$ complexity.
The methods we describe in this section can approximate the same dense linear algebra calculation in $\mathcal{O}(N)$
or $\mathcal{O}(N\log N)$ time.

As an example of the absolute performance of the analytic variants,
we refer to the Treecode implementation -- \texttt{Bonsai},
which scales to the full node of Titan using 18,600 GPUs achieving 24.77 PFlops \cite{Bedorf2014}.
Bonsai's performance comes not only from its matrix-free nature,
but also from domain specific optimizations
for hardcoded quadrupoles and an assumption that all charges are positive.
Therefore, this kind of performance cannot be transferred to other applications that require higher accuracy.
However, viewing these methods as a preconditioner instead of a direct solver
significantly reduces the accuracy requirements \cite{Ibeid2016,Aminfar2016a}.

\subsection{Fast Translation Operators}
A large part of the calculation time of FMM is spent on the translation
of multipole expansions to local expansions (or their equivalent charges).
Therefore, much work has focused on developing fast translation operators to accelerate this part of the FMM.
Rotation of spherical harmonics \cite{White1996}, Block FFT \cite{Elliott1996}, Planewaves \cite{Greengard1997}
are analytic options for fast translation operators.

These translation operators are applied to a pair of boxes
in the FMM tree structure that satisfy a certain proximity threshold.
This proximity is usually defined as the parent's neighbors' children that are non-neighbors.
This produces a list of boxes that are far enough that the multipole/local expansion converges,
but are close enough that the expansion does not converge for the their parents.
Such an interaction list can contain up to $6^3-3^3=189$ source boxes for each target box.
Out of these 189 boxes, the ones that are further from the target box can perform
the translation operation using their parent box as the source without loss of accuracy.
There are a few variants for these techniques that reduce the interaction list size
such as the level-skip M2L method \cite{Wang2015}
and 8,4,2-box method \cite{Wilkes2015}.
There are also methods that use the dual tree traversal along with the multipole acceptance criterion
to construct optimal interaction lists \cite{Dehnen2002},
which automates the process of finding the optimal interaction list size.

Another technique to accelerate the translation operators is the use of variable expansion order,
as proposed in the very fast multipole method (VFMM) \cite{Petersen1994}, Gaussian VFMM \cite{Burant1996},
optimal parameter FMM \cite{Choi2001}, and error controlled FMM \cite{Dachsel2010}.
There are two main reasons why spatially varying the expansion order in the translation operators is beneficial.
One is because not all boxes in the interaction list are of equal distance,
and the boxes that are further from each other can afford to use lower expansion order,
while retaining the accuracy.
The other reason is because some parts of the domain may have smaller values,
and the contribution from that part can afford to use lower expansion order
without sacrificing the overall accuracy.

The translation operators can be stored as matrices that operate on the vector of expansion coefficients.
Therefore, singular value decomposition (SVD) can be used to compress this matrix \cite{Gimbutas2002}
and BLAS can be used to maximize the cache utilization \cite{Fortin2005}.
Some methods use a combination of these techniques like Chebychev with SVD \cite{Fong2009}
and planewave with adaptive cross approximation (ACA) and SVD \cite{Hesford2011}.
The use of SVD is a systematic and optimal way of achieving what the variable expansion order techniques
in the previous paragraph were trying to do manually.
Precomputing these translation matrices and storing them is a typical optimization technique
in many FMM implementations \cite{Malhotra2015}.

One important connection to make here is that these matrices for the translation operators are precisely
what $\mathcal{H}^2$-matrices and $HSS$ matrices store in the off-diagonal blocks after compression.
One can think of FMM as a method that has the analytical form to generate these small matrices
in the off-diagonal blocks, without relying on numerical low-rank approximation methods.
To complete this analogy, we point out that the dense diagonal blocks in $\mathcal{H}^2$-matrices and $HSS$ matrices
are simply storing the direct operator (Green's function) in FMM.
Noticing this equivalence leads to many possibilities of hybridization among the analytic and algebraic variants.
Possibly the most profound is the following.
Those that are familiar with FMM know that translation operators for boxes with the same relative positioning are identical.
This suggests that many of the entries in the off-diagonal blocks of
$\mathcal{H}^2$-matrices and $HSS$ matrices are identical.
For matrices that are generated from a mesh that has a regular structure even the diagonal blocks would be identical,
which is what happens in FMMs for continuous volume integrals \cite{Malhotra2015}.
This leads to $\mathcal{O}(1)$ storage for the matrix entries at every level of the hierarchy,
so the total storage cost of these hierarchical matrices could be reduced to $\mathcal{O}(\log N)$
if the identical entires are not stored redundantly.
This aspect is currently underutilized in the algebraic variants, but seems obvious from the analytic side.
By making use of the translational invariance and rotational symmetry of the interaction list
one can reduce the amount of storage even further \cite{Coulaud2008,Darve2011,Takahashi2012}.
This also results in blocking techniques for better cache utilization.

\subsection{Semi-analytical FMM}
The methods described in the previous subsection all require the existance of an analytical form
of the multipole/local translation operator, which is kernel dependent.
There are a class of methods that remove this restriction by using equivalent charges
instead of multipole expansions \cite{Anderson1992,Berman1995,Makino1999}.
A well known implementation of this method is the kernel independent FMM (KIFMM) code \cite{Ying2004}.
There are also variants that use Chebychev polynomials \cite{Dutt1996},
and a representative implementation of this is the Black-box FMM \cite{Fong2009}.
As the name of these codes suggest, these variants of the FMM have reduced requirements for the
information that has to be provided by the user.
The translation operators are kernel-independent, which frees the user from the most difficult task
of having to provide an analytical form of the translation operators.
For example, if one wants to calculate the Mat\'{e}rn function for covariance martices,
or multiquadrics for radial basis function interpolation,
one simply needs to provide these functions and the location of the points and the FMM will handle the rest.
It is important to note that these methods are not entirely kernel independent or black-box
because the user still needs to provide the kernel dependent analytic form of the
original equation they wish to calculate.
Using the vocabulary of the algebraic variants, one could say that these analytical expressions for
the hierarchical matrices are kernel independent only for the off-diagonal blocks,
and for the diagonal blocks the analytical form is kernel dependent.

FMM for continuous volume integrals \cite{Ethridge2001} also has important features
when considering the analytic-algebraic tradeoff.
The volume integrals are often combined with boundary integrals, as well \cite{Ying2006}.
One can think of these methods as an FMM that includes the discretization process \cite{Langston2011}.
Unlike the FMM for discrete particles, these methods have the ability to impose regular underlying geometry.
This enables the use of precomputation of the direct interaction matrix in the analytic variants \cite{Malhotra2015},
and reduces the storage requirements of the dense diagonal blocks in the algebraic variants.

\begin{table}[b]
\caption{Categorization of algebraic low-rank approximation methods.}
\begin{tabular}{p{3cm}p{2.5cm}p{3cm}p{2cm}}
\hline\noalign{\smallskip}
Method & Hierarchical & Weak admissibility & Nested basis \\
\noalign{\smallskip}\svhline\noalign{\smallskip}
$\mathcal{H}$-matrix \cite{Hackbusch1999} & yes & maybe & no \\
$\mathcal{H}^2$-matrix \cite{Hackbusch2000a} & yes & maybe & yes \\
HODLR \cite{Ambikasaran2013} & yes & yes & no \\
HSS \cite{Chandrasekaran2006} / HBS \cite{Martinsson2005} & yes & yes & yes \\
BLR \cite{Amestoy2015} & no & yes & no \\
\noalign{\smallskip}\hline\noalign{\smallskip}
\end{tabular}
\label{tab:algebraic}
\end{table}

\subsection{Algebraic Low-Rank Approximation}
There are many variants of algebraic low-rank approximation methods.
They can be categorized based on whether they are hierarchical, whether they use weak admissibility,
or if the basis is nested, as shown in Table \ref{tab:algebraic}.
For the definition of admissibility see \cite{Grasedyck2003}.
Starting from the top, $\mathcal{H}$-matrices \cite{Hackbusch1999,Bebendorf2008} are hierarchical,
usually use standard or strong admissibility, and no nested basis.
The analytic counterpart of the $\mathcal{H}$-matrix is the Treecode.
The $\mathcal{H}^2$-matrices \cite{Hackbusch2000a,Borm2009} are also hierarchical
and use standard or strong admissibility, but unlike $\mathcal{H}$-matrices use a nested basis.
This brings the complexity down from $\mathcal{O}(NlogN)$ to $\mathcal{O}(N)$.
The analytic counterpart of the $\mathcal{H}^2$-matrix is the FMM.
The next three entries in Table \ref{tab:algebraic} do not have analytic counterparts because
analytic low-rank approximations do not converge under weak admissibility conditions.
Hierarchical off-diagonal low-rank (HODLR) matrices \cite{Ambikasaran2013,Aminfar2016},
are basically $\mathcal{H}$-matrices with weak admissibility conditions.
Similarly, hierarchically semi-seperable (HSS) \cite{Chandrasekaran2006, Xia2010},
and hierarchically block-seperable (HBS) \cite{Martinsson2005} matrices
are $\mathcal{H}^2$-matrices with weak admissibility conditions.
The block low-rank (BLR) matrices \cite{Amestoy2015} are a non-hierarchical version of the HODLR,
with just the bottom level.
A summary of implementations and their characteristics are presented in \cite{Rouet2015}.

For methods that do not have weak admissibility, it is common to use geometrical information to calculate
the standard/strong admissibility condition.
This dependence on the geometry of the algebraic variants is not ideal.
There have been various proposals for algebraic clustering methods \cite{LeBorne2006,Oliveira2007,Grasedyck2008}.
This problem requires even more advanced solutions for high dimension problems \cite{March2015a}.
Stronger admissibility is also problem for parallelization since it results in more communication.
There have been studies on how to partition hierarchical matrices on distributed memory \cite{Izadi2012}.
There are also methods to reduce the amount of memory consumption
during the construction of HSS matrices\cite{Lessel2015}.  
  
The categorization in Table \ref{tab:algebraic} is for the hierarchical matrix structure,
and any low-rank approximation method can be used with each of them during the compression phase.
The singular value decomposition is the most na\"{i}ve and expensive way to calculate
a low-rank approximation.
QR or LU decompositions can be used to find the numerical rank by using appropriate pivoting.
Rank-revealing QR \cite{Chan1987} has been proposed along with efficient pivoting strategies
\cite{Hong1992, Chandrasekaran1994, Gu1996}.
Rank-revealing LU \cite{Chan1984} also requires efficient pivoting strategies \cite{Hwang1992,Hwang1997,Miranian2003}.
Rank-revealing LU is typically faster than rank-revealing QR \cite{Pan2000}.
There are other methods like the pseudo-skeletal method \cite{Goreinov1997}
and adaptive cross approximation (ACA) \cite{Bebendorf2000,Bebendorf2003},
which do not yield the optimal low-rank factorizations but have a much lower cost.
ACA has a better pivoting strategy than pseudo-skeletal methods,
but can still fail because of bad pivots \cite{Borm2003}.
The hybrid cross approximation (HCA) \cite{Borm2005} has the same proven convergence
as standard interpolation but also the same efficiency as ACA.
Yet another class of low-rank approximation is the interpolative decomposition (ID) \cite{Cheng2005a, Martinsson2005},
where a few of its columns are used to form a well-conditioned basis for the remaining columns.
ID can be combined with randomized methods \cite{Liberty2007}, which has much lower complexity.
For a nice review on these randomized methods see \cite{Halko2011}.

\section{Low-Rank Approximation for Factorization}
\subsection{Sparse Matrix Factorization}
Hierarchical low-rank approximation methods can be used as direct solvers with controllable accuracy.
This makes them useful as preconditioners within a Krylov subspace method,
which in turn reduces the accuracy requirements of the low-rank approximation.
High accuracy and completely algebraic methods are demanding in terms of memory consumption
and amount of communication, so they will not be the optimal choice if they are not the only option for that problem.
It is worth noting that these methods are for solving dense matrices,
and for sparse matrices it is best to combine them with multifrontal methods
and use them to compress the Schur complements \cite{Xia2009}.
They should never be used to compress a sparse matrix nor the inverse of it directly,
even if the inverse of a sparse matrix is dense.
There are various methods to reduce fill-in during the factorization to not make the inverse dense,
and not using these methods first will have a toll on the asymptotic constant \cite{Grasedyck2009},
even though the asymptotic complexity may still be optimal.

Ultimately, minimizing fill-in and minimizing off-diagonal rank should not be conflicting objectives.
The former depends on the connectivity and the latter depends on the distance in the underlying geometry.
In most applications, the closer points are connected (or interact) more densely,
so reordering according to the distance should produce near optimal ordering for the connectivity as well.
The same can be said about minimizing communication for the parallel implementation of these methods.
Mapping the 3-D connectivity/distance to a 1-D locality in the memory space (or matrix column/row)
is what we are ultimately trying to achieve.

There are various ways to minimize the fill-in and compress the dense blocks during factorization.
These dense blocks (Schur complements) are an algebraic form of the Green's function \cite{Xia2014},
and have the same low-rank properties \cite{Chandrasekaran2010} stemming from the fact that
some of the boundary points in the underlying geometry are distant from each other.
Formulating a boundary integral equation is the analytical way of arriving to the same dense matrix.
From an algebraic point of view, the sparse matrix for the volume turns into a dense matrix for the boundary,
through the process of trying to minimize fill-in.
Considering the minimization of fill-in and the compression of the dense matrices in separate phases leads to
methods like HSS + multifrontal \cite{Xia2009,Xia2010,Xia2013}.

\subsection{Dense Matrix Factorization}
The methods in the previous subsection are direct solvers/preconditioners for sparse matrices.
As we have mentioned, there is an analogy between minimizing fill-in in sparse matrices by looking at the connectivity,
and minimizing the rank of off-diagonal blocks of dense matrices by looking at the distance.
Using this analogy, the same concept as nested dissection for sparse matrices can be applied to dense matrices.
This leads to methods like the recursive skeletonization \cite{Ho2012},
or hierarchical Poincare-Steklov (HPS) \cite{Martinsson2015,Gillman2015}.
HPS is like a bottom-up version of what nested dissection and recursive skeletonization do top-down.
For high contrast coefficient problems, it makes sense to construct the domain dissection bottom-up,
to align the bisectors with the coefficient jumps.
There are also other methods that rely on a similar concept \cite{Yarvin1999,Greengard2009,Kong2011,Bremer2012a}.
Furthermore, since many of these methods use weak admissibility with growing ranks for 3-D problems,
it is useful to have nested hierarchical decompositions, which is like a nested dimension reduction.
In this respect, the recursive skeletonization has been extended to
hierarchical interpolative factorization (HIF) \cite{Ho2015},
the HSS has been extended to HSS2D \cite{Xia2014}.
There is also a combination of HSS and Skeletonization \cite{Corona2015}.
There are methods that use this nested dimension reduction concept without the low-rank approximation \cite{Henon2006}
in the context of domain decomposition for incomplete LU factorization.
One method that does not use weak admissibility is the inverse FMM \cite{Ambikasaran2014},
which makes it applicable to 3-D problems in $\mathcal{O}(N)$ without nested dimension reduction.

\section{Experimental Results}
\subsection{FMM vs. HSS}
There have been very few comparisons between the analytic and algebraic hierarchical low-rank approximation methods.
From a high performance computing perspective, the practical performance of highly optimized implementations
of these various methods is of great interest.
There have been many efforts to develop new methods in this area, which has resulted in a large amount
of similar methods with different names without a clear overall picture of their relative performance
on modern HPC architectures.
The trend in architecture where arithmetic operations are becoming cheap compared to data movement,
is something that must be considered carefully when predicting which method will perform better on computers of the future.

We acknowledge that the comparisons we present here are far from complete,
and much more comparisons between all the different methods are needed in order to acheive our long term objective.
The limitation actually comes from the lack of highly optimized implementations of these methods
that are openly available to us at the moment.

In the present work we start by comparing exaFMM -- a highly optimized implementation of FMM,
with STRUMPACK -- a highly optimized implementation of HSS.
We select the 2D and 3D Laplace equation on uniform lattices as test cases.
For HSS we directly construct the compressed matrix by calling the Green's function in
the randomized low-rank approximation routine.
We perform the matrix-vector multiplication using the FMM and HSS,
and measure the time for the compression/precalculation and application of the matrix-vector multiplication.
We also measure the peak memory consumption of both methods.

\begin{figure}[t]
\centering
\includegraphics[scale=.55]{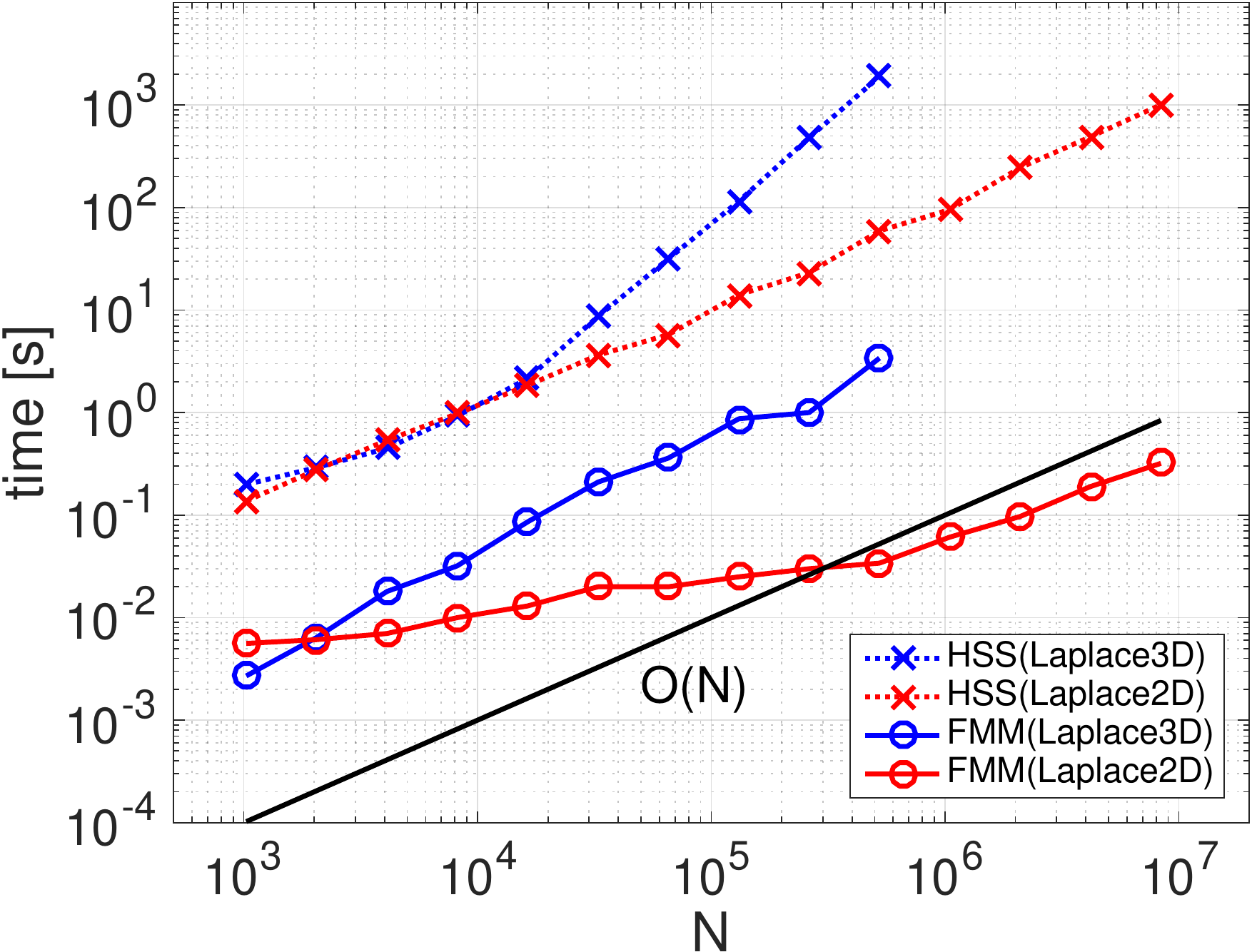}
\caption{Elapsed time for the matrix-vector multiplication using FMM and HSS for different problem sizes.}
\label{fig:fmm_hss}
\end{figure}

The elapsed time for the FMM and HSS for different problem sizes is shown in Fig. \ref{fig:fmm_hss}.
In order to isolate the effect of the thread scalability of the two methods,
these runs are performed on a single core of a 12-core Ivy Bridge (E5-2695 v2).
For the 2D Laplace equation, the FMM shows some overhead for small $N$,
but is about 3 orders of magnitude faster than HSS for larger problems.
For the 3D Laplace equation, the FMM is about 2 orders of magnitude faster than HSS for smaller $N$,
but HSS exhibits non-optimal behavior for large $N$ because the rank keeps growing.

\begin{figure}[t]
\centering
\includegraphics[scale=.5]{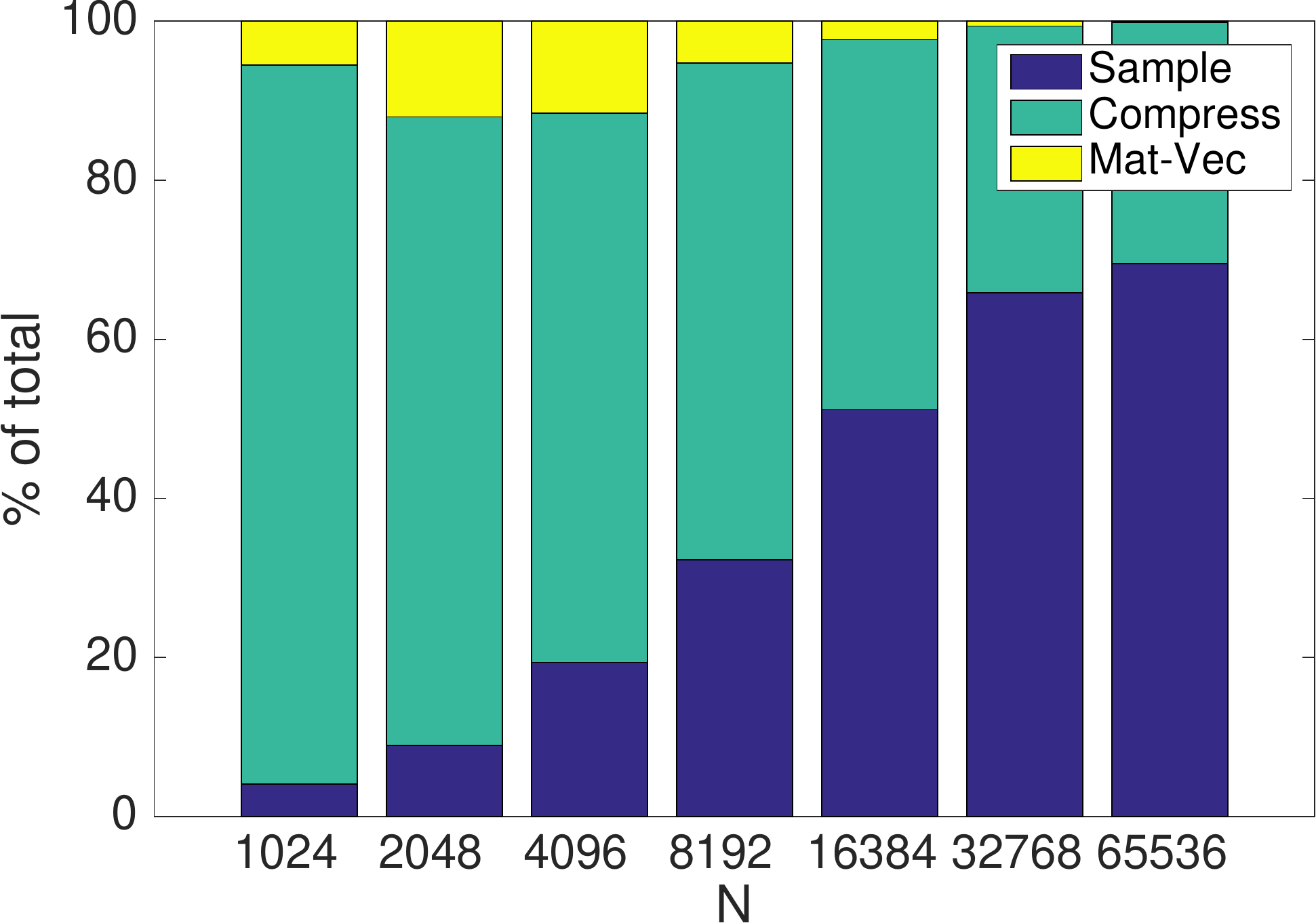}
\caption{Percentage of the computation time of HSS for different problem sizes.}
\label{fig:breakdown}
\end{figure}

The large difference in the computational time is actually coming from the heavy computation in the
sampling phase and compression phase of the HSS.
In Fig. \ref{fig:breakdown}, we show the percentage of the computation time of HSS for different problem sizes $N$.
``Sample" is the sampling time, ``Compress" is the compression time, and ``Mat-Vec" is the matrix-vector multiplication time.
We can see that the sampling is taking longer and longer as the problem size increases.
This is because the rank $k$ increases with the problem size $N$,
and both sampling and compression time increase with the $k$ and $N$.

\begin{figure}[h]
\centering
\includegraphics[scale=.55]{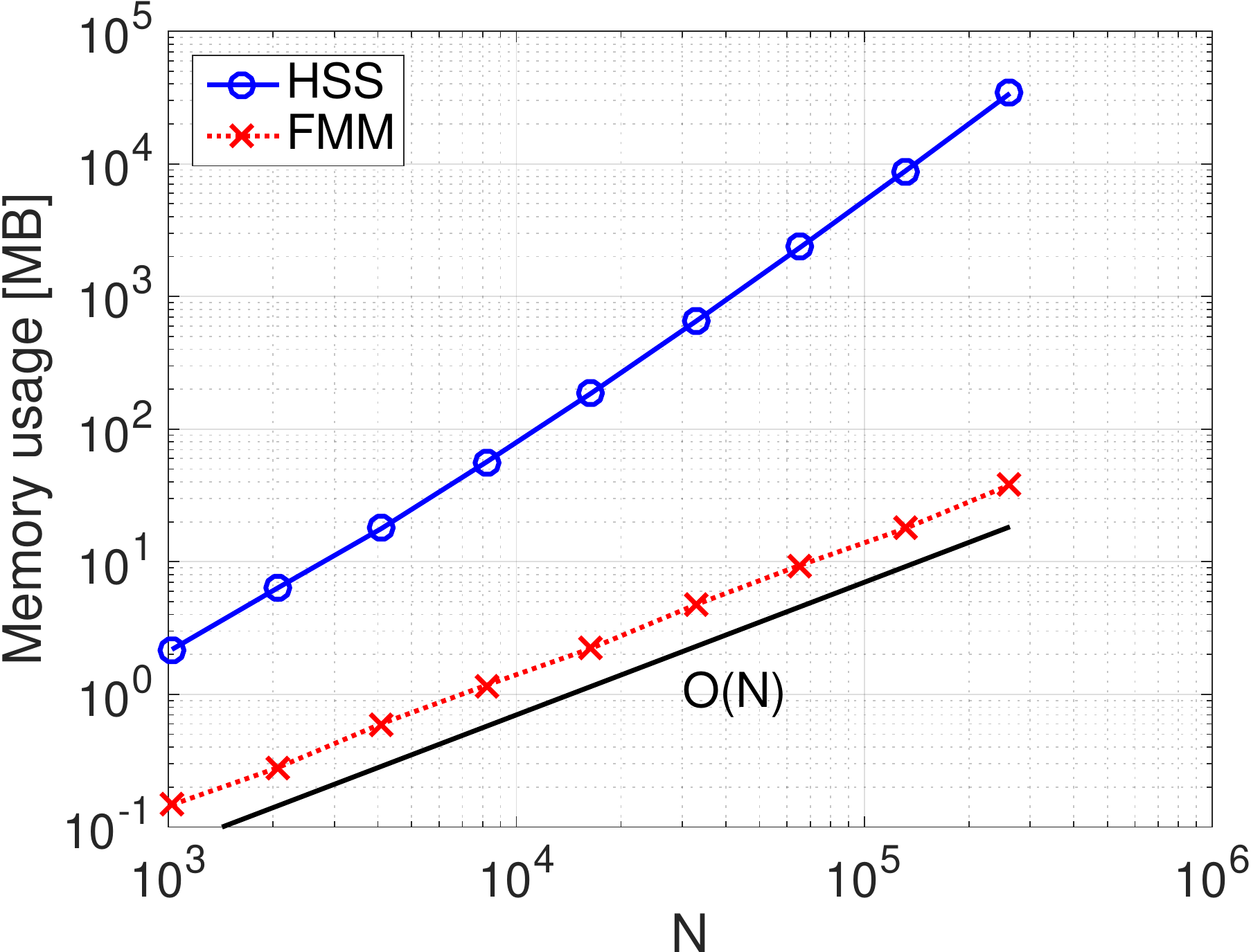}
\caption{Peak memory usage of FMM and HSS for the 3D Laplace equation.}
\label{fig:memory}
\end{figure}

The peak memory usage of FMM and HSS is shown in Fig. \ref{fig:memory} for the 3D Laplace equation.
We see that the FMM has strictly $\mathcal{O}(N)$ storage requirements,
but since the rank in the HSS grows for 3D kernels it does not show the ideal $\mathcal{O}(N\log N)$ behavior.
The disadvantage of HSS is two-fold.
First of all, its algebraic nature requires it to store the compressed matrix,
where as the FMM is analytic and therefore matrix-free.
Secondly, the weak admissibility causes the rank to grow for 3D problems,
and with that the memory consumption grows at a suboptimal complexity.

\subsection{FMM vs. Multigrid}
If we are to use the FMM as a matrix-free $\mathcal{O}(N)$ preconditioner based on hierarchical low-rank approximation,
the natural question to ask is "How does it compare against multigrid?",
which is a much more popular matrix-free $\mathcal{O}(N)$ preconditioner for solving elliptic PDEs.
We perform a benchmark test similar to the one in the previous subsection,
for the Laplace equation and Helmholtz equation on a 3D cubic lattice $[-1,1]^3$,
but for this case we impose Dirichlet boundary conditions at the faces of the domain.
The preconditioners are used inside a Krylov subspace solver.

\begin{figure}[t]
\centering
\includegraphics[scale=.7]{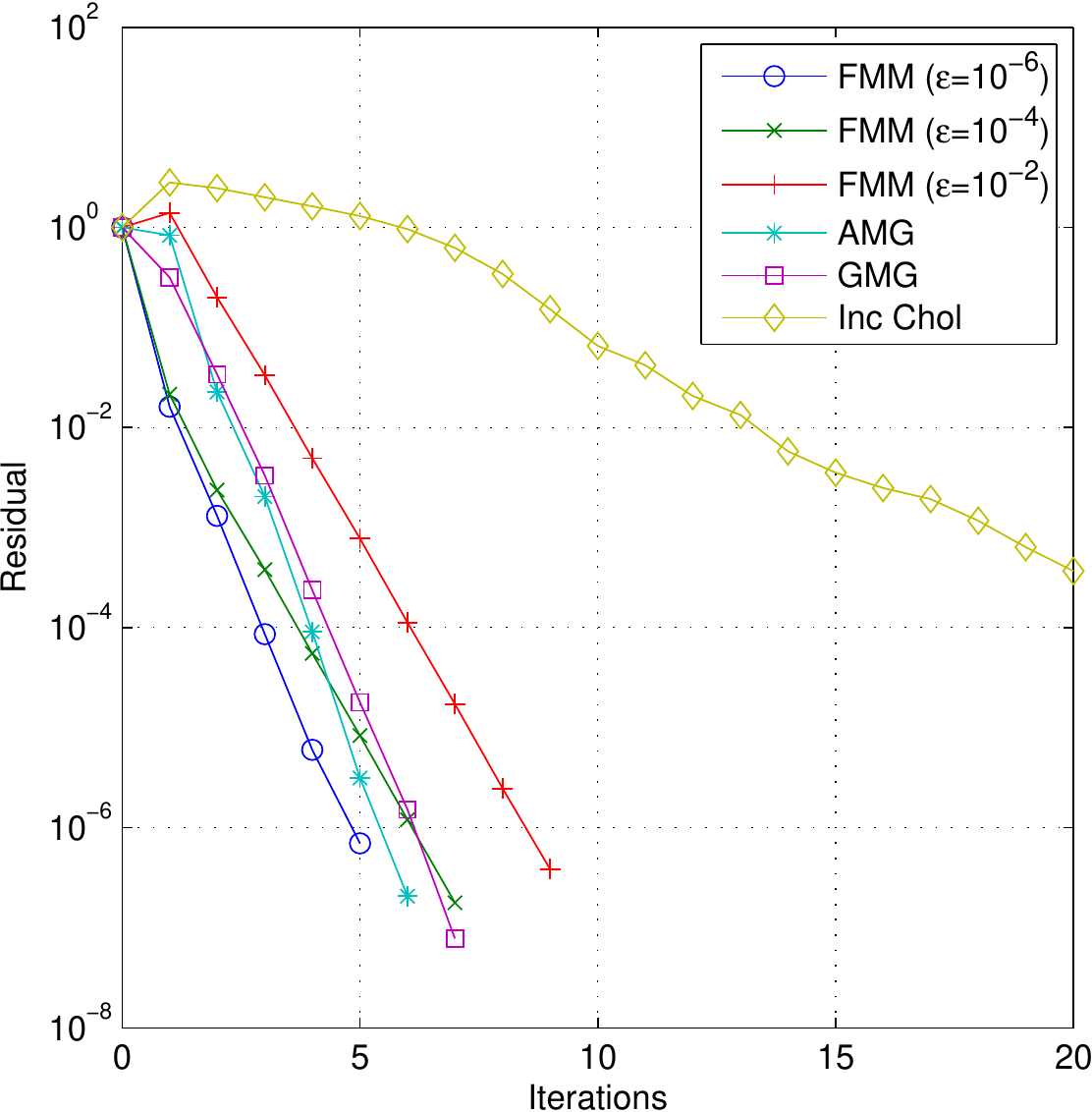}
\caption{Convergence rate of the FMM and Multigrid preconditioners for
the Laplace equation on a $[-1,1]^3$ lattice with spacing $h=2^{-5}$.}
\label{fig:laplace}
\end{figure}

\begin{figure}[t]
\centering
\includegraphics[scale=.7]{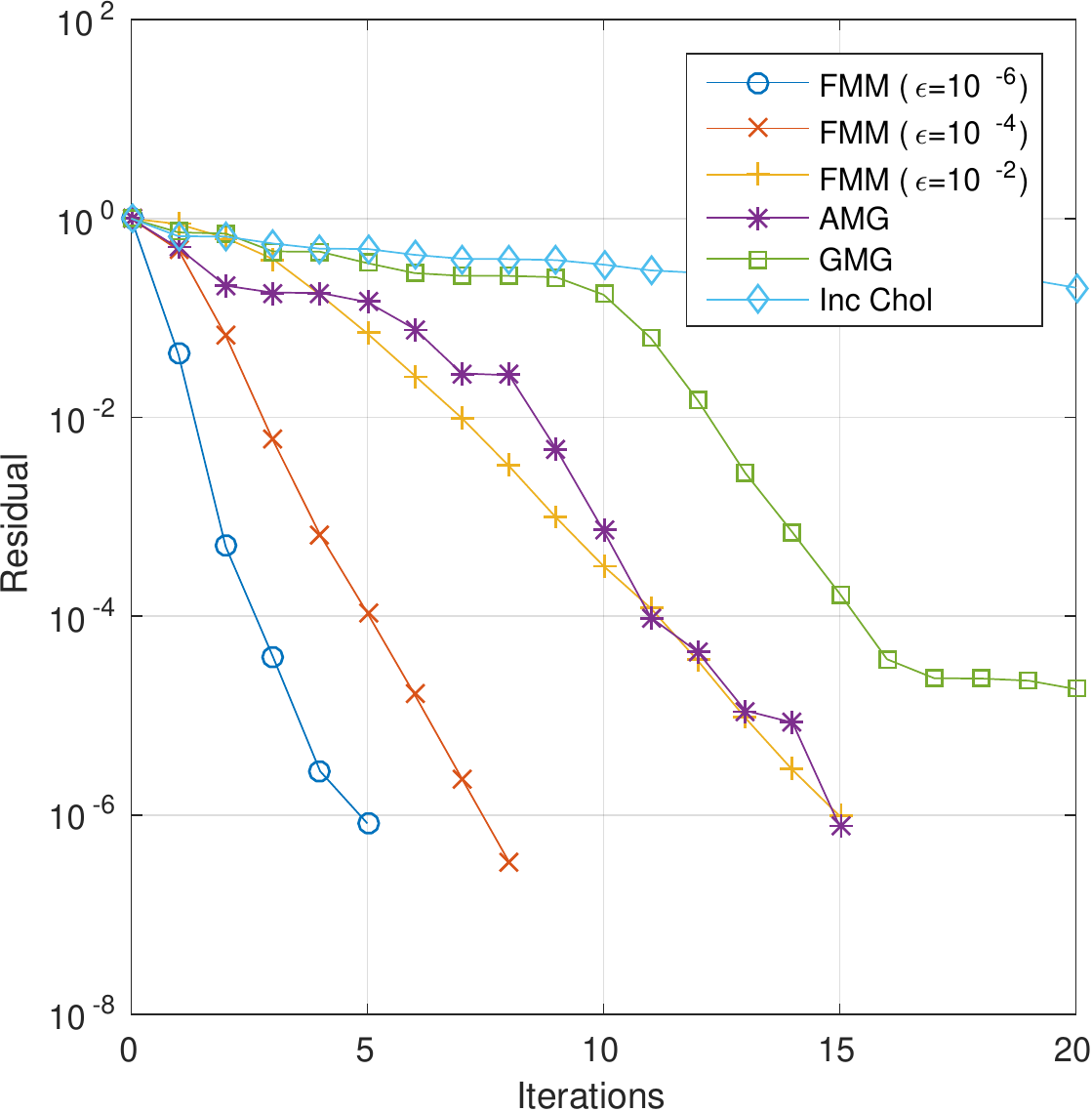}
\caption{Convergence rate of the FMM and Multigrid preconditioners for the Helmholtz equation on a $[-1,1]^3$ lattice with
spacing $h=2^{-5}$ and wave number $\kappa=7$.}
\label{fig:helmholtz}
\end{figure}

The convergence rate of the FMM and Multigrid preconditioners for the Laplace equation is shown in Fig. \ref{fig:laplace},
for a grid spacing of $h=2^{-5}$.
``AMG" is algebraic multigrid, ``GMG" is geometric multigrid, ``Inc Chol" is incomplete Cholesky.
The $\epsilon$ value represents the accuracy of the FMM.
We see that the FMM preconditioner has comparable convergence to the algebraic and geometric multigrid method.
Even for a very low-accuracy FMM with $\epsilon=10^{-2}$,
the convergence rate is much better than the incomplete Cholesky.
We refer to the work by Ibeid \textit{et al.} \cite{Ibeid2016} for more detailed comparisons between FMM and Multigrid.

A similar plot is shown for the Helmholtz equation with grid spacing of $h=2^{-5}$
and wave number $\kappa=7$ in Fig. \ref{fig:helmholtz}.
The nomenclature of the legend is identical to that of Fig. \ref{fig:laplace}.
In this case, we see a larger difference between the convergence rate of FMM and Multigrid.
Even the FMM with the worst accuracy does better than the multigrid.
We have also confirmed that the FMM preconditioner has a convergence rate that is independent of the problem size,
up to moderate wave numbers of $\kappa$.

\section{Conclusions and Outlook}
We have shown the contrast between the analytical and algebraic hierarchical low-rank approximations,
by reviewing the contributions over the years and placing them along the analytical-algebraic spectrum.
The relation between Treecode, FMM, KIFMM, black-box FMM,
$\mathcal{H}$-matrix, $\mathcal{H}^2$-matrix, HODLR, HSS,
HBS, and BLR were explained from the perspective of compute-memory tradeoff.
This birds-eye view of the entire hierarchical low-rank approximation landscape from analytical to algebraic,
allows us to place ideas like precomputation of FMM translation matrices
and relate that to storage reduction techniques for the algebraic variants.

Some important findings from this cross-disciplinary literature review are:
\begin{itemize}
\item Translational invariance of the FMM operators suggest that $\mathcal{H}^2$-matrices (and the like)
have mostly duplicate entries, which many are redundantly storing at the moment.
\item The analytical variants can now perform factorization and are kernel independent,
so the decision to use the algebraic variants at the cost of consuming more memory
should be made carefully.
\item The kernel-independent variants of FMM can be used as a matrix-free $\mathcal{O}(N)$ compression technique.
\item The use of SVD to compress the FMM translation matrices,
makes the work on variable expansion order and its error optimized variants redundant.
\item The hierarchical compression should not be applied directly to the inverse or factorizations of sparse matrices
just because they fill-in. One must first try to minimize fill-in, and then compress only the dense blocks that cannot be avoided.
\end{itemize}

The comparison benchmarks between FMM and HSS are still preliminary tests for a very simple case.
However, they clearly demonstrate the magnitude of the difference that lies between the various
hierarchical low-rank approximation methods.
The comparison between FMM and multigrid is also a very simple test case,
but it reveals the previously unquantified convergence properties of low-accuracy FMM as a preconditioner.
Of course, for such simple problems the FMM can give the exact solution in
finite arithmetic and therefore solve the problem in a single iteration.
The interesting point here is not the fact that it can be used as a preconditioner,
but the practical performance of the low-accuracy FMM being significantly faster than the high accuracy FMM,
even if it requires a few iterations.

There is much more that can be done if all of these complicated hierarchical low-rank approximation methods
could somehow be made easier to code.
We believe a modular view of these methods will help the developers though separation of concerns.
Instead of everyone coding a slightly different version of the whole thing,
we could each choose a module to focus on that fits our research interests,
and contribute to a larger and more sustainable ecosystem.
A few ideas to facilitate the transition to such a community effort are:
\begin{enumerate}
\item Create a common benchmark (mini app) for each of the modules.
\item Gradually propagate standards in the community, starting from the major codes.
\item Develop a common interface between the hierarchical structure and inner kernels.
\item Do not try to unify code, just have a standard with a common API (like MPI).
\end{enumerate}

\begin{acknowledgement}
We thank Fran\c{c}ois-Henry Rouet, Pieter Ghysels, and Xiaoye, S. Li for providing the STRUMPACK interface for our comparisons between FMM and HSS. This work was supported by JSPS KAKENHI Grant-in-Aid for Research Activity Start-up Grant Number 15H06196. This publication was based on work supported in part by Award No KUK-C1-013-04, made by King Abdullah University of Science and Technology (KAUST). This work used the Extreme Science and Engineering Discovery Environment (XSEDE), which is supported by National Science Foundation grant number OCI-1053575. 
\end{acknowledgement}

\bibliographystyle{abbrv}

\begin{thebibliography}{100}

\bibitem{Ambikasaran2013}
S.~Ambikasaran and E.~Darve.
\newblock An {O(NlogN)} fast direct solver for partial hierarchically
  semi-separable matrices.
\newblock {\em Journal of Scientific Computing}, 57:477--501, 2013.

\bibitem{Ambikasaran2014}
S.~Ambikasaran and E.~Darve.
\newblock The inverse fast multipole method.
\newblock {\em arXiv:1407.1572v1}, 2014.

\bibitem{Ambikasaran2013a}
S.~Ambikasaran, J.-Y. Li, P.~K. Kitanidis, and E.~Darve.
\newblock Large-scale stochastic linear inversion using hierarchical matrices.
\newblock {\em Computational Geosciences}, 17(6):913--927, 2013.

\bibitem{Amestoy2015}
P.~Amestoy, C.~Ashcraft, O.~Boiteau, A.~Buttari, J.-Y. L'Excellent, and
  C.~Weisbecker.
\newblock Improving multifrontal methods by means of block low-rank
  representations.
\newblock {\em SIAM Journal on Scientific Computing}, 37(3):A1451--A1474, 2015.

\bibitem{Aminfar2016}
A.~Aminfar, S.~Ambikasaran, and E.~Darve.
\newblock A fast block low-rank dense solver with applications to
  finite-element matrices.
\newblock {\em Journal of Computational Physics}, 304:170--188, 2016.

\bibitem{Aminfar2016a}
A.~Aminfar and E.~Darve.
\newblock A fast, memory efficient and robust sparse preconditioner based on a
  multifrontal approach with applications to finite-element matrices.
\newblock {\em International Journal for Numerical Methods in Engineering},
  accepted, 2016.

\bibitem{Anderson1992}
C.~R. Anderson.
\newblock An implementation of the fast multipole method without multipoles.
\newblock {\em SIAM Journal on Scientific and Statistical Computing},
  13(4):923--947, 1992.

\bibitem{Appel1985}
A.~W. Appel.
\newblock An efficient program for many-body simulation.
\newblock {\em SIAM Journal on Scientific and Statistical Computing},
  6(1):85--103, 1985.

\bibitem{Barba2013}
L.~A. Barba and R.~Yokota.
\newblock How will the fast multipole method fare in the exascale era?
\newblock {\em SIAM News}, 46(6):1--3, 2013.

\bibitem{Barnes1986}
J.~Barnes and P.~Hut.
\newblock {O(NlogN)} force-calculation algorithm.
\newblock {\em Nature}, 324:446--449, 1986.

\bibitem{Bebendorf2000}
M.~Bebendorf.
\newblock Approximation of boundary element matrices.
\newblock {\em Numerische Mathematik}, 86:565--589, 2000.

\bibitem{Bebendorf2008}
M.~Bebendorf.
\newblock {\em Hierarchical Matrices}, volume~63 of {\em Lecture Notes in
  Computational Science and Engineering}.
\newblock Springer, 2008.

\bibitem{Bebendorf2003}
M.~Bebendorf and S.~Rjasanow.
\newblock Adaptive low-rank approximation of collocation matrices.
\newblock {\em Computing}, 70:1--24, 2003.

\bibitem{Bedorf2014}
J.~B{\'e}dorf, E.~Gaburov, M.~S. Fujii, K.~Nitadori, T.~Ishiyama, and
  S.~Portegies~Zwart.
\newblock 24.77 {P}flops on a gravitational tree-code to simulate the milky way
  galaxy with 18600 {GPU}s.
\newblock In {\em Proceedings of the 2014 ACM/IEEE International Conference for
  High Performance Computing, Networking, Storage and Analysis}, pages 1--12,
  2014.

\bibitem{Berman1995}
C.~L. Berman.
\newblock Grid-multipole calculations.
\newblock {\em SIAM Journal on Scientific Computing}, 16(5):1082--1091, 1995.

\bibitem{Borm2009}
S.~B{\"o}rm.
\newblock Construction of data-sparse $h^2$-matrices by hierarchical
  compression.
\newblock {\em SIAM Journal on Scientific Computing}, 31(3):1820--1839, 2009.

\bibitem{Borm2005}
S.~B{\"o}rm and L.~Grasedyck.
\newblock Hybrid cross approximation of integral operators.
\newblock {\em Numerische Mathematik}, 101:221--249, 2005.

\bibitem{Borm2003}
S.~B{\"o}rm, L.~Grasedyck, and W.~Hackbusch.
\newblock Introduction to hierarchical matrices with applications.
\newblock {\em Engineering Analysis with Boundary Elements}, 27:405--422, 2003.

\bibitem{Bremer2012a}
J.~Bremer.
\newblock A fast direct solver for the integral equations of scattering theory
  on planar curves with corners.
\newblock {\em Journal of Computational Physics}, 231:1879--1899, 2012.

\bibitem{Burant1996}
J.~C. Burant, M.~C. Strain, G.~E. Scuseria, and M.~J. Frisch.
\newblock Analytic energy gradients for the {G}aussian very fast multipole
  method ({GvFMM}).
\newblock {\em Chemical Physics Letters}, 248:43--49, 1996.

\bibitem{Chaillat2008}
S.~Chaillat, M.~Bonnet, and J.-F. Semblat.
\newblock A multi-level fast multipole {BEM} for 3-{D} elastodynamics in the
  frequency domain.
\newblock {\em Computer Methods in Applied Mechanics and Engineering},
  197:4233--4249, 2008.

\bibitem{Chan1984}
T.~F. Chan.
\newblock On the existence and computation of {LU}-factorizations with small
  pivots.
\newblock {\em Mathematics of Computation}, 42(166):535--547, 1984.

\bibitem{Chan1987}
T.~F. Chan.
\newblock Rank revealing {QR} factorizations.
\newblock {\em Linear Algebra and its Applications}, 88/89:67--82, 1987.

\bibitem{Chandrasekaran2006}
S.~Chandrasekaran, P.~Dewilde, M.~Gu, W.~Lyons, and T.~Pals.
\newblock A fast solver for {HSS} representations via sparse matrices.
\newblock {\em SIAM Journal on Matrix Analysis and Applications}, 29(1):67--81,
  2006.

\bibitem{Chandrasekaran2010}
S.~Chandrasekaran, P.~Dewilde, M.~Gu, and N.~Somasunderam.
\newblock On the numerical rank of the off-diagonal blocks of {S}chur
  complements of discretized elliptic {PDE}s.
\newblock {\em SIAM Journal on Matrix Analysis and Applications},
  31(5):2261--2290, 2010.

\bibitem{Chandrasekaran1994}
S.~Chandrasekaran and I.~C.~F. Ipsen.
\newblock On rank-revealing factorizations.
\newblock {\em SIAM Journal on Matrix Analysis and Applications},
  15(2):592--622, 1994.

\bibitem{Cheng2005a}
H.~Cheng, Z.~Gimbutas, P.~G. Martinsson, and V.~Rokhlin.
\newblock On the compression of low rank matrices.
\newblock {\em SIAM Journal on Scientific Computing}, 26(4):1389--1404, 2005.

\bibitem{Choi2001}
C.~H. Choi, K.~Ruedenberg, and M.~S. Gordon.
\newblock New parallel optimal-parameter fast multipole method ({OPFMM}).
\newblock {\em Journal of Computational Chemistry}, 22(13):1484--1501, 2001.

\bibitem{Corona2015}
E.~Corona, P.~G. Martinsson, and D.~Zorin.
\newblock An {O(N)} direct solver for integral equations on the plane.
\newblock {\em Applied and Computational Harmonic Analysis}, 38:284--317, 2015.

\bibitem{Coulaud2008}
O.~Coulaud, P.~Fortin, and J.~Roman.
\newblock High performance {BLAS} formulation of the multipole-to-local
  operator in the fast multipole method.
\newblock {\em Journal of Computational Physics}, 227:1836--1862, 2008.

\bibitem{Dachsel2010}
H.~Dachsel.
\newblock Corrected article: ``an error-controlled fast multipole method".
\newblock {\em The Journal of Chemical Physics}, 132:119901, 2010.

\bibitem{Darve2011}
E.~Darve, C.~Cecka, and T.~Takahashi.
\newblock The fast multipole method on parallel clusters, multicore processors,
  and graphics processing units.
\newblock {\em Comptes Rendus Mecanique}, 339:185--193, 2011.

\bibitem{Darve2004a}
E.~Darve and P.~Hav{\'e}.
\newblock A fast multipole method for {M}axwell equations stable at all
  frequencies.
\newblock {\em Philosophical Transactions of the Royal Society of London A},
  362:603--628, 2004.

\bibitem{Dehnen2002}
W.~Dehnen.
\newblock A hierarchical {O(N)} force calculation algorithm.
\newblock {\em Journal of Computational Physics}, 179(1):27--42, 2002.

\bibitem{Dutt1996}
A.~Dutt, M.~Gu, and V.~Rokhlin.
\newblock Fast algorithms for polynomial interpolation, integration, and
  differntiation.
\newblock {\em SIAM Journal on Numerical Analysis}, 33(5):1689--1711, 1996.

\bibitem{Elliott1996}
W.~D. Elliott and J.~A. Board.
\newblock Fast {F}ourier transform accelerated fast multipole algorithm.
\newblock {\em SIAM Journal on Scientific Computing}, 17(2):398--415, 1996.

\bibitem{Ethridge2001}
F.~Ethridge and L.~Greengard.
\newblock A new fast-multipole accelerated {P}oisson solver in two dimensions.
\newblock {\em SIAM Journal on Scientific Computing}, 23(3):741--760, 2001.

\bibitem{Fong2009}
W.~Fong and E.~Darve.
\newblock The black-box fast multipole method.
\newblock {\em Journal of Computational Physics}, 228:8712--8725, 2009.

\bibitem{Fortin2005}
P.~Fortin.
\newblock Multipole-to-local operator in the fast multipole method:
  {C}omparison of {FFT}, rotations and {BLAS} improvements.
\newblock Technical Report RR-5752, Rapports de recherche, et theses de
  l'Inria, 2005.

\bibitem{Gillman2015}
A.~Gillman, A.~Barnett, and P.~G. Martinsson.
\newblock A spectrally accurate direct solution technique for frequency-domain
  scattering problems with variable media.
\newblock {\em BIT Numerical Mathematics}, 55:141--170, 2015.

\bibitem{Gimbutas2013}
Z.~Gimbutas and L.~Greengard.
\newblock Fast multi-particle scattering: A hybrid solver for the {M}axwell
  equations in microstructured materials.
\newblock {\em Journal of Computational Physics}, 232:22--32, 2013.

\bibitem{Gimbutas2002}
Z.~Gimbutas and V.~Rokhlin.
\newblock A generalized fast multipole method for nonoscillatory kernels.
\newblock {\em SIAM Journal on Scientific Computing}, 24(3):796--817, 2002.

\bibitem{Goreinov1997}
S.~A. Goreinov, E.~E. Tyrtyshnikov, and N.~L. Zamarashkin.
\newblock A theory of pseudoskeleton approximations.
\newblock {\em Linear Algebra and its Applications}, 261(1-3):1--21, 1997.

\bibitem{Grasedyck2003}
L.~Grasedyck and W.~Hackbusch.
\newblock Construction and arithmetics of {H}-matrices.
\newblock {\em Computing}, 70:295--334, 2003.

\bibitem{Grasedyck2008}
L.~Grasedyck, R.~Kriemann, and S.~Le~Borne.
\newblock Parallel black box {H-LU} preconditioning for elliptic boundary value
  poblems.
\newblock {\em Computing and Visualization in Science}, 11:273--291, 2008.

\bibitem{Grasedyck2009}
L.~Grasedyck, R.~Kriemann, and S.~Le~Borne.
\newblock Domain decomposition based {H-LU} preconditioning.
\newblock {\em Numerische Mathematik}, 112:565--600, 2009.

\bibitem{Gray2001}
A.~G. Gray and A.~W. Moore.
\newblock {N}-body problems in statistical learning.
\newblock In T.~K. Leen, T.~G. Dietterich, and V.~Tresp, editors, {\em Advances
  in Neural Information Processing Systems}, volume~13, pages 521---527. MIT
  Press, 2001.

\bibitem{Greengard2009}
L.~Greengard, D.~Gueyffier, P.~G. Martinsson, and V.~Rokhlin.
\newblock Fast direct solvers for integral equations in complex three
  dimensional domains.
\newblock {\em Acta Numerica}, 18:243--275, 2009.

\bibitem{Greengard1987}
L.~Greengard and V.~Rokhlin.
\newblock A fast algorithm for particle simulations.
\newblock {\em Journal of Computational Physics}, 73(2):325--348, 1987.

\bibitem{Greengard1997}
L.~Greengard and V.~Rokhlin.
\newblock A new version of the fast multipole method for the {L}aplace equation
  in three dimensions.
\newblock {\em Acta Numerica}, 6:229--269, 1997.

\bibitem{Gu1996}
M.~Gu and E.~S. C.
\newblock Efficient algorithms for computing a strong rank-revealing {QR}
  factorization.
\newblock {\em SIAM Journal on Scientific Computing}, 17(4):848--869, 1996.

\bibitem{Gumerov2007}
N.~A. Gumerov and R.~Duraiswami.
\newblock Fast radial basis function interpolation via preconditioned {K}rylov
  iteration.
\newblock {\em SIAM Journal on Scientific Computing}, 29(5):1876--1899, 2007.

\bibitem{Hackbusch1999}
W.~Hackbusch.
\newblock A sparse matrix arithmetic based on {H}-matrices, part {I}:
  {I}ntroduction to {H}-matrices.
\newblock {\em Computing}, 62:89--108, 1999.

\bibitem{Hackbusch2000a}
W.~Hackbusch, B.~Khoromskij, and S.~A. Sauter.
\newblock On $h^2$-matrices.
\newblock In H.~Bungartz, R.~Hoppe, and C.~Zenger, editors, {\em Lectures on
  Applied Mathematics}. Springer-Verlag, 2000.

\bibitem{Hackbusch1989}
W.~Hackbusch and Z.~P. Nowak.
\newblock On the fast matrix multiplication in the boundary element method by
  panel clustering.
\newblock {\em Numerische Mathematik}, 54:463--491, 1989.

\bibitem{Halko2011}
N.~Halko, P.~G. Martinsson, and J.~A. Tropp.
\newblock Finding structure with randomness: Probabilistic algorithms for
  constructing approximate matrix decompositions.
\newblock {\em SIAM Review}, 53(2):217--288, 2011.

\bibitem{Hao2015}
S.~Hao, P.~G. Martinsson, and P.~Young.
\newblock An efficient and highly accurate solver for multi-body acoustic
  scattering problems involving rotationally symmetric scatterers.
\newblock {\em Computers and Mathematics with Applications}, 69:304--318, 2015.

\bibitem{Henon2006}
P.~H{\'e}non and Y.~Saad.
\newblock A parallel multistage {ILU} factorization based on a hierarchical
  graph decomposition.
\newblock {\em SIAM Journal on Scientific Computing}, 28(6):2266--2293, 2006.

\bibitem{Hesford2011}
A.~J. Hesford and R.~C. Waag.
\newblock Reduced-rank approximations to the far-field transform in the gridded
  fast multipole method.
\newblock {\em Journal of Computational Physics}, 230:3656--3667, 2011.

\bibitem{Ho2012}
K.~L. Ho and L.~Greengard.
\newblock A fast direct solver for structured linear systems by recursive
  skeletonization.
\newblock {\em SIAM Journal on Scientific Computing}, 34(5):A2507--A2532, 2012.

\bibitem{Ho2015}
K.~L. Ho and L.~Ying.
\newblock Hierarchical interpolative factorization for elliptic operators:
  Integral equations.
\newblock {\em arXiv:1307.2666}, 2015.

\bibitem{Hong1992}
Y.~P. Hong and C.~T. Pan.
\newblock Rank-revealing {QR} factorizations and the singular value
  decomposition.
\newblock {\em Mathematics of Computation}, 58(197):213--232, 1992.

\bibitem{Hwang1997}
T.-M. Hwang, W.-W. Lin, and D.~Pierce.
\newblock Improved bound for rank revealing {LU} factorizations.
\newblock {\em Linear Algebra and its Applications}, 261(1):173--186, 1997.

\bibitem{Hwang1992}
T.-M. Hwang, W.-W. Lin, and E.~K. Yang.
\newblock Rank revealing {LU} factorizations.
\newblock {\em Linear Algebra and its Applications}, 175:115--141, 1992.

\bibitem{Ibeid2016}
H.~Ibeid, R.~Yokota, J.~Pestana, and D.~Keyes.
\newblock Fast multipole preconditioners for sparse matrices arising from
  elliptic equations.
\newblock {\em arXiv:1308.3339}, 2016.

\bibitem{Izadi2012}
M.~Izadi.
\newblock {\em Hierarchical Matrix Techniques on Massively Parallel Computers}.
\newblock PhD thesis, Universitat Leipzig, 2012.

\bibitem{Kong2011}
W.~Y. Kong, J.~Bremer, and V.~Rokhlin.
\newblock An adaptive fast direct solver for boundary integral equations in two
  dimensions.
\newblock {\em Applied and Computational Harmonic Analysis}, 31:346--369, 2011.

\bibitem{Langston2011}
H.~Langston, L.~Greengard, and D.~Zorin.
\newblock A free-space adaptive {FMM}-based {PDE} solver in three dimensions.
\newblock {\em Communications in Applied Mathematics and Computational
  Science}, 6(1):79--122, 2011.

\bibitem{LeBorne2006}
S.~Le~Borne.
\newblock Multilevel hierarchical matrices.
\newblock {\em SIAM Journal on Matrix Analysis and Applications},
  28(3):871--889, 2006.

\bibitem{Lee2012}
D.~Lee, R.~Vuduc, and A.~G. Gray.
\newblock A distributed kernel summation framework for general-dimension
  machine learning.
\newblock In {\em Proceedings of the 2012 SIAM International Conference on Data
  Mining}, 2012.

\bibitem{Lessel2015}
K.~Lessel, M.~Hartman, and S.~Chandrasekaran.
\newblock A fast memory efficient construction algorithm for hierarchically
  semi-separable representations.
\newblock {\em http://scg.ece.ucsb.edu/publications/MemoryEfficientHSS.pdf},
  2015.

\bibitem{Li2014}
J.-Y. Li, S.~Ambikasaran, E.~F. Darve, and P.~K. Kitanidis.
\newblock A {K}alman filter powered by $h^2$-matrices for quasi-continuous data
  assimilation problems.
\newblock {\em Water Resources Research}, 50:3734--3749, 2014.


\bibitem{Liang2013}
Z.~Liang, Z.~Gimbutas, L.~Greengard, J.~Huang, and S.~Jiang.
\newblock A fast multipole method for the {R}otne-{P}rager-{Y}amakawa tensor
  and its applications.
\newblock {\em Journal of Computational Physics}, 234:133--139, 2013.

\bibitem{Liberty2007}
E.~Liberty, F.~Woolfe, P.~G. Martinsson, V.~Rokhlin, and M.~Tygert.
\newblock Randomized algorithms for the low-rank approximation of matrices.
\newblock {\em PNAS}, 104(51):20167--20172, 2007.

\bibitem{Makino1999}
J.~Makino.
\newblock Yet another fast multipole method without multipoles --
  {P}seudoparticle multipole method.
\newblock {\em Journal of Computational Physics}, 151(2):910--920, 1999.

\bibitem{Malhotra2015}
D.~Malhotra and G.~Biros.
\newblock {PVFMM}: A parallel kernel independent {FMM} for particle and volume
  potentials.
\newblock {\em Communications in Computational Physics}, 18(3):808--830, 2015.

\bibitem{Malhotra2014}
D.~Malhotra, A.~Gholami, and G.~Biros.
\newblock A volume integral equation stokes solver for problems with variable
  coefficients.
\newblock In {\em Proceedings of the 2014 ACM/IEEE International Conference for
  High Performance Computing, Networking, Storage and Analysis}, pages 1--11,
  2014.

\bibitem{March2015a}
W.~B. March, B.~Xiao, and G.~Biros.
\newblock {ASKIT}: Approximate skeletonization kernel-independent treecode in
  high dimensions.
\newblock {\em SIAM Journal on Scientific Computing}, 37(2):A1089--A1110, 2015.

\bibitem{Martinsson2015}
P.~G. Martinsson.
\newblock The hierarchical {P}oincar{\'e}-{S}teklov ({HPS}) solver for elliptic
  {PDE}s: A tutorial.
\newblock {\em arXiv:1506.01308}, 2015.

\bibitem{Martinsson2005}
P.~G. Martinsson and V.~Rokhlin.
\newblock A fast direct solver for boundary integral equations in two
  dimensions.
\newblock {\em Journal of Computational Physics}, 205:1--23, 2005.

\bibitem{Miranian2003}
L.~Miranian and M.~Gu.
\newblock Strong rank revealing {LU} factorizations.
\newblock {\em Linear Algebra and its Applications}, 367:1--16, 2003.

\bibitem{Ohno2014}
Y.~Ohno, R.~Yokota, H.~Koyama, G.~Morimoto, A.~Hasegawa, G.~Masumoto,
  N.~Okimoto, Y.~Hirano, H.~Ibeid, T.~Narumi, and M.~Taiji.
\newblock Petascale molecular dynamics simulation using the fast multipole
  method on k computer.
\newblock {\em Computer Physics Communications}, 185:2575--2585, 2014.

\bibitem{Oliveira2007}
S.~Oliveira and Y.~F.
\newblock An algebraic approcah for {H}-matrix preconditioners.
\newblock {\em Computing}, 80:169--188, 2007.

\bibitem{Pan2000}
C.~T. Pan.
\newblock On the existence and computation of rank-revealing {LU}
  factorizations.
\newblock {\em Linear Algebra and its Applications}, 316:199--222, 2000.

\bibitem{Petersen1994}
H.~G. Petersen, D.~Soelvason, J.~W. Perram, and E.~R. Smith.
\newblock The very fast multipole method.
\newblock {\em The Journal of Chemical Physics}, 101(10):8870--8876, 1994.

\bibitem{Rahimian2010}
A.~Rahimian, I.~Lashuk, K.~Veerapaneni, A.~Chandramowlishwaran, D.~Malhotra,
  L.~Moon, R.~Sampath, A.~Shringarpure, J.~Vetter, R.~Vuduc, D.~Zorin, and
  G.~Biros.
\newblock Petascale direct numerical simulation of blood flow on 200k cores and
  heterogeneous architectures.
\newblock In {\em Proceedings of the 2010 ACM/IEEE International Conference for
  High Performance Computing, Networking, Storage and Analysis}, SC '10, 2010.

\bibitem{Rouet2015}
F.-H. Rouet, X.-S. Li, P.~Ghysels, and A.~Napov.
\newblock A distributed-memory package for dense hierarchically semi-separable
  matrix computations using randomization.
\newblock {\em arXiv:1503.05464}, 2015.

\bibitem{Shao2001}
Y.~Shao, C.~A. White, and M.~Head-Gordon.
\newblock Efficient evaluation of the {C}oulomb force in density-functional
  theory calculations.
\newblock {\em The Journal of Chemical Physics}, 114(15):6572--6577, 2001.

\bibitem{Takahashi2012}
T.~Takahashi, C.~Cecka, W.~Fong, and E.~Darve.
\newblock Optimizing the multipole-to-local operator in the fast multipole
  method for graphical processing units.
\newblock {\em International Journal for Numerical Methods in Engineering},
  89:105--133, 2012.

\bibitem{Verde2015}
A.~Verde and A.~Ghassemi.
\newblock Fast multipole displacement discontinuity method ({FM-DDM}) for
  geomechanics reservoir simulations.
\newblock {\em International Journal for Numerical and Analytical Methods in
  Geomechanics}, 39(18):1953--1974, 2015.

\bibitem{Wang2015}
Y.~Wang, Q.~Wang, X.~Deng, Z.~Xia, J.~Yan, and H.~Xu.
\newblock Graphics processing unit ({GPU}) accelerated fast multipole {BEM}
  with level-skip {M2L} for {3D} elasticity problems.
\newblock {\em Advances in Engineering Software}, 82:105--118, 2015.

\bibitem{White1996}
C.~A. White and M.~Head-Gordon.
\newblock Rotating around the quartic angular momentum barrier in fast
  multipole method calculations.
\newblock {\em The Journal of Chemical Physics}, 105(12):5061--5067, 1996.

\bibitem{Wilkes2015}
D.~R. Wilkes and A.~J. Duncan.
\newblock A low frequency elastodynamic fast multipole boundary element method
  in three dimensions.
\newblock {\em Computational Mechanics}, 56:829--848, 2015.

\bibitem{Willis2005}
D.~Willis, J.~Peraire, and J.~White.
\newblock {FastAero} -- a precorrected {FFT}-fast multipole tree steady and
  unsteady potential flow solver.
\newblock {\em http://hdl.handle.net/1721.1/7378}, 2005.

\bibitem{Wolf2011a}
W.~R. Wolf and S.~K. Lele.
\newblock Aeroacoustic integrals accelerated by fast multipole method.
\newblock {\em AIAA Journal}, 49(7):1466--1477, 2011.

\bibitem{Xia2013}
J.~Xia.
\newblock Randomized sparse direct solvers.
\newblock {\em SIAM Journal on Matrix Analysis and Applications},
  34(1):197--227, 2013.

\bibitem{Xia2014}
J.~Xia.
\newblock {O(N)} complexity randomized 3{D} direct solver with {HSS2D}
  structure.
\newblock Proceedings of the Project Review, Geo-Mathematical Imaging Group
  317--325, Purdue University, 2014.

\bibitem{Xia2009}
J.~Xia, S.~Chandrasekaran, M.~Gu, and X.~S. Li.
\newblock Superfast multifrontal method for large structured linear systems of
  equations.
\newblock {\em SIAM Journal on Matrix Analysis and Applications},
  31(3):1382--1411, 2009.

\bibitem{Xia2010}
J.~Xia, S.~Chandrasekaran, M.~Gu, and X.~S. Li.
\newblock Fast algorithms for hierarchically semiseperable matrices.
\newblock {\em Numerical Linear Algebra with Applications}, 17:953--976, 2010.

\bibitem{Yarvin1999}
N.~Yarvin and V.~Rokhlin.
\newblock An improved fast multipole algorithm for potential fields on the
  line.
\newblock {\em SIAM Journal on Numerical Analysis}, 36(2):629--666, 1999.

\bibitem{Ying2004}
L.~Ying, G.~Biros, and D.~Zorin.
\newblock A kernel-independent adaptive fast multipole algorithm in two and
  three dimensions.
\newblock {\em Journal of Computational Physics}, 196(2):591--626, 2004.

\bibitem{Ying2006}
L.~Ying, G.~Biros, and D.~Zorin.
\newblock A high-order 3{D} boundary integral equation solver for elliptic
  {PDE}s in smooth domains.
\newblock {\em Journal of Computational Physics}, 219:247--275, 2006.

\bibitem{Yokota2011}
R.~Yokota, J.~P. Bardhan, M.~G. Knepley, L.~A. Barba, and T.~Hamada.
\newblock Biomolecular electrostatics using a fast multipole {BEM} on up to 512
  {GPU}s and a billion unknowns.
\newblock {\em Computer Physics Communications}, 182:1272--1283, 2011.

\bibitem{Yokota2013}
R.~Yokota, T.~Narumi, K.~Yasuoka, and L.~A. Barba.
\newblock Petascale turbulence simulation using a highly parallel fast
  multipole method on {GPU}s.
\newblock {\em Computer Physics Communications}, 184:445--455, 2013.

\bibitem{Yunis2012}
E.~Yunis, R.~Yokota, and A.~Ahmadia.
\newblock Scalable force directed graph layout algorithms using fast multipole
  methods.
\newblock In {\em The 11th International Symposium on Parallel and Distributed
  Computing}, Munich, Germany, June 2012.

\bibitem{Zhao2007}
Z.~Zhao, N.~Kovvali, W.~Lin, C.-H. Ahn, L.~Couchman, and L.~Carin.
\newblock Volumetric fast multipole method for modeling {S}chr{\"o}dinger's
  equation.
\newblock {\em Journal of Computational Physics}, 224:941--955, 2007.

\end{thebibliography}

\end{document}